\documentclass[a4paper,leqno]{article}

\usepackage{amsmath,amssymb,amscd,cite,authblk,multicol,color}

\allowdisplaybreaks

\begin{document}

\begin{center}

\textbf{\Large Projective (or spin) representations of finite groups. II}

\bigskip

{Tatsuya Tsurii\footnote{Tokyo University of Information Sciences, 
Onaridai 4-1, Wakaba-ku, Chiba-shi, Chiba, 265-8501, JAPAN}, 
Satoe Yamanaka\footnote{Department of Liberal Studies, 
National Institute of Technology, Nara College, 22 Yata-cho, Yamatokoriyama, Nara 639-1080, JAPAN}, 
Itsumi Mikami\footnote{Hirai Mathematics Institute}, and 
Takeshi Hirai\footnote{Kyoto University, Yoshida-honmachi, Sakyo-ku, Kyoto 606-8501, JAPAN}
\footnote{Hirai Mathematics Institute, 22-5 Nakazaichi-cho, Iwakura, Sakyo, Kyoto 606-0027, JAPAN} }

\bigskip

\end{center}

\begin{abstract}
In the previous paper, we proposed a practical method of 
constructing 
explicitly 
representation groups $R(G)$ for finite groups $G$, 
and apply it 
to certain  typical finite groups $G$ with 
Schur multiplier $M(G)$ containing prime number 3. 
In this paper, we construct a complete list of irreducible 
projective (or spin) representations of $G$ and compute their 
characters (called spin characters). 
It is a continuation of our study  of spin representations 
in the cases 
where $M(G)$ contains prime number 2  to 
the cases where other prime $p$ appears, firstly $p=3$. 
We classify  
irreducible spin representations and calculate 
spin characters according 
to their spin types.

\bigskip

\medskip


\end{abstract}

%%%%%%%%%%%%%%%%%%%%%%%%%%%%%%%%%%%%%%%%%%%%%%%%%%%%%%%%%%%%%%55 sect 1
\begin{multicols*}{2}

{\bf 1. Introduction. } 
Let $G$ be a finite group. Representation group 
  $R(G)$ is a certain central extension of $G$ by $M(G)$ as its 
central subgroup, and any  
irreducible spin representation $\pi$ of $G$ can be 
induced naturally from 
a linear representation $\widetilde{\pi}$ of $R(G)$.   
There exists  a one-dimensional character 
$\chi$ of $M(G)$ such that 
$\widetilde{\pi}(z)=\chi(z)I\;\;(z\in M(G))$, called  
 {\it spin type} of $\pi$, where $I$ denotes the identity operator. 

Our practical 
method of constructing representation group $R(G)$ 
of $G$, 
named as \lq\lq{\it Efficient stairway up 
to the Sky}\,\rq\rq, is 
to start with a presentation of $G$ by a pair of {\bf Set of generators} 
and {\bf Set of fundamental relations}. A central extension 
\begin{gather}
\label{2023-11-29-1}
\quad 
1 \longrightarrow A \longrightarrow H
\stackrel{\alpha}{\longrightarrow} G\longrightarrow 1\quad ({\rm exact})
\end{gather}
of $G$ by an abelian group $A$ is called 
an {\it efficient covering} if 
 $A$ is contained in 
the commutator subgroup $[H,H]$ of $H$, i.e., 
\begin{gather}
\label{2023-11-29-2}
A\subset Z(H)\cap [H,H].
\end{gather}
Take a pair $x,y\in G$ commuting with each other and of the same 
order $d>1$. Then try to take a {\it one-step efficient 
covering} given by the 
commuting relation $[x,y]=1$, in such a way that,  $A=\langle z\rangle$ 
is a cyclic group generated by 
 a central element $z$ of order $d$ which is given under the homomorphism 
$\alpha$ as   
\begin{gather}
\label{2023-11-29-3}
[\xi,\eta]=z,\quad \xi\stackrel{\alpha}{\rightarrow}  x,\;
\eta \stackrel{\alpha}{\rightarrow} y,\;z\stackrel{\alpha}{\rightarrow}1.
\end{gather}

When we know Schur multiplier $M(G)$ of $G$, we may arrive up to $R(G)$ 
by repeating one-step efficient central extensions step by step because  
$R(G)$ is characterized \cite[\S 2]{HMTY} 
as an efficient central extension $H$ 
of $G$ such that
$$
M(G) \subset Z(H)\cap [H,H], \quad H/M(G)\cong G.
$$

In cases we are treating here, for 
every one-step efficient central extensions, 
the resulting central extension 
$H$ is expressed as semidirect product $H=U\rtimes W$, 
and 
 so, construction 
of its irreducible representations is realized by Mackey's induced 
 representations when $U$ is abelian 
(cf. \cite[Theorem 5.1]{YTMH}), 
and by repeated application 
of the classical method 
 in \cite{THirai} for general semidirect products.

{\bf 2. Preliminaries.}
Let $G'$ be a central extension of $G$ by an abelian group  $Z$ as 
$$
1 \rightarrow Z \rightarrow G'
 \stackrel{\alpha}{\rightarrow} G \rightarrow 1.
$$
Take a section $s:G \to {\cal S}_G\subset G'$ 
for the canonical homomorphism $\alpha$.  Then, 
 we have a $Z$-valued cocycle 
\,$z_{g,h}\in Z\;(g,h\in G)$\,  such that 
$$
s(g)s(h)=z_{g,h}s(gh),\quad \exists \,z_{g,h}\in Z.
$$ 
For an irreducible linear representation  
$\Pi$ of $G'$, put\; 
$$
\pi(g):= \Pi\big(s(g)\big)\quad(g\in G). 
$$ 
Let 
$\chi\in\widehat{Z}$ be the {\it spin type}\, of 
$\Pi$ with respect to $Z$, that is, 
\,$
\Pi(z)=\chi(z)I\;\,(z\in Z)
$, 
then $\pi$ is a spin representation of $G$ associated with a 
factor set 
$$
r_{g,h}= \chi(z_{g,h})\in 
{\boldsymbol T}^1:=\{\lambda\in 
{\boldsymbol C}; |\lambda|=1\}.
$$ 
In fact, for $g, h\in G$, we have 
\begin{gather*}
\pi(g)\pi(h)=\Pi\big(s(g)\big)\Pi\big(s(h)\big)
=
\Pi\big(s(g)s(h)\big)
\\
=\Pi(z_{g,h}s(gh)\big)=\chi(z_{g,h})\Pi(s(gh)\big)
=r_{g,h}\pi(gh).
\end{gather*}
We call $\pi$ a {\it sectional restriction} (onto $G$) of $\Pi$.

If we take another section $s': G \to G'$, 
then it induces a spin representation 
$\pi'(g):=\Pi\big(s'(g)\big)$ of $G$ 
associated to an equivalent factor set 
$r'_{g,h}\;(g,h\in G)$. 

We know by Schur the following. 

\vskip.8ex

{\bf Theorem 2.1.} (Schur \cite{Sch1})\; {\it 
As a central extension $G'$ above, 
take a representation 
group $R(G)$. Then, for any factor set 
$r_{g,h}$ for $G$, there exists an irreducible 
linear representation\, $\Pi$ of $R(G)$ which 
induces as above spin representation 
$\pi$ of\, $G$ associated to a factor set which 
is equivalent to the given $r_{g,h}$. 
}%\it  

\vskip.8ex

{\bf 3. Irreducible representations of semidirect product groups. }  
Let $G$ be a compact group of semidirect product type\;  
$
G=U\rtimes W,
$ 
where $U$ is a compact group, normal in $G$, and $W$ 
is a finite group.  We gave in \cite{THirai} 
 the so-called {\it classical method}\,  
to construct a complete list of representatives of the 
dual $\widehat{G}$, or of 
all the equivalence classes of 
IRs of $G$. We  
 will apply it in the case of finite groups. 

So, 
starting from now on, we assume $G$ to be finite, and explain 
briefly about the method.  
\vskip.8ex

{\bf 3.1. Construction of IRs of 
 semidirect product groups.} 
Take an IR $\rho$ of $U$ and consider its equivalence 
class $[\rho]\in\widehat{U}$. 
Every $w\in W$ acts on $u\in U$ as $u\to w(u)$, and on $\rho$  
as 
$$
({}^w\!\rho)(u):=\rho\big(w^{-1}(u)\big)\quad(u\in U),
$$ 
and therefore 
 on equivalent classes as $[\rho]\to [{}^w\!\rho]$. 
Denote by $\widehat{U}/W$ the set of $W$-orbits 
in the dual $\widehat{U}$ of $U$. 

{\bf (Step 1)} 
Take a stationary subgroup $W([\rho])$ in $W$ of an 
equivalence class $[\rho]$, 
that is, 
$$
W([\rho]):=\{w\in W\,;\,
{}^w\!\rho\cong \rho\},
$$
 and put $H:=U\rtimes W([\rho])$. 

{\bf (Step 2)}  
For $w \in W([\rho])$, we determine 
explicitly an intertwining operator $J_\rho(w)$ given as   
\begin{gather}
\label{2023-11-30-1}
\quad
\rho\big(w(u)\big)=J_\rho(w)\,\rho(u)\,J_\rho(w)^{-1}\;\;(u\in U). 
\end{gather}
Then it is determined up to a non-zero scalar factor.  
 Hence we have a spin representation 
$W([\rho])\ni w \mapsto J_\rho(w)$. Let $\alpha_{v,w}$ 
be its factor set given as \;
 $$
 J_\rho(v)J_\rho(w)= \alpha_{v,w}\,J_\rho(vw)\quad 
 \big(v,w \in W([\rho])\big).
 $$

{\bf (Step 3)} 
Let $W([\rho])'$ be a central extension of
 $W([\rho])$ sufficiently large compared to 
the cocycle $\alpha_{v,w}\,$:\;   
$$
\quad
1\to Z \to W([\rho])' \stackrel{\Phi}{\longrightarrow} 
W([\rho]) \to 1
\quad{\rm  (exact)},
$$ 
where $\Phi$ denotes the canonical homomorphism. 
Then, by Theorem 2.1, $J_\rho$ can be lifted up to 
a linear representation $J'_\rho$ of $W([\rho])'$. 
  Put $H':=U\rtimes W([\rho])'$ with the action 
$$
\;\;w'(u):=w(u),\;\;w'\in W([\rho])',\;w=\Phi(w'). 
$$
Put also 
\begin{eqnarray}
\label{2023-11-30-2}
\nonumber
~\pi^0\big((u,w')\big):=\rho(u)\cdot J'_\rho(w')
\;\;\big(u\in U,\,w'\in W([\rho])'\big), 
\end{eqnarray}
then $\pi^0=\rho\cdot J'_\rho$ is an IR of $H'$.

{\bf (Step 4)}  
Take an IR $\pi^1$ of $W([\rho])'$ 
and consider it as a representation of $H'$ 
through the homomorphism 
\;$
H' \to\, W([\rho])' \cong H'/U, 
$\; 
and consider inner tensor product 
$$
\pi:=\pi^0\boxdot \pi^1
$$ 
as a representation of $H'$. 
Let the factor set of $\pi^1$, 
viewed as a spin representation of the base group $W([\rho])$, 
be $\beta_{v,w}$\,, then that of 
$\pi$ is 
\,$
\alpha_{v,w}\,\beta_{v,w}. 
$

 To get an IR of $G$, we should pick up $\pi^1$ with 
the factor set $\beta_{v,w}=\alpha_{v,w}^{\;-1}$. 
This is possible by Theorem 2.1. Then $\pi$ 
becomes a linear representation of the base group 
$H=U\rtimes W([\rho])$.  
 Thus we obtain a representation of $G$ by inducing it up as 
\begin{gather}
\label{2023-11-30-11}
\quad
\Pi(\pi^0,\pi^1) :={\rm Ind}^G_H\,\pi={\rm Ind}^G_H(\pi^0\boxdot \pi^1). 
\end{gather}

{\bf Theorem 3.1.}\; (\hspace{-0.1mm}\cite[Theorem 3.1]{THirai})
{\it 
 Let $G=U\rtimes W$ be as above a semidirect product 
group. 
Then the induced representation\; 
$
\Pi(\pi^0,\pi^1)={\rm Ind}^G_H (\pi^0\boxdot\pi^1)
$\;  
of\, $G$ in\, {\rm (\ref{2023-11-30-11})}\, is irreducible. 

}%\it  

\vskip.8ex
 
{\bf 3.2. Character of $\Pi(\pi^0,\pi^1)$.} 
 Put $\Pi=\Pi(\pi^0,\pi^1)$, and let $\chi_\Pi$ be the trace 
character of $\Pi$. Since $\Pi={\rm Ind}_H^G\pi$, it is  
 expressed as 
\begin{eqnarray*}
\label{2023-11-30-3}
\chi_\Pi(g) =\int_{H\backslash G}\chi_\pi(kgk^{-1})\,d\nu_{H\backslash G}(\dot{k}) 
\end{eqnarray*}
where the character $\chi_\pi$ of $\pi$ is extended 
from $H$ to $G$ by putting 0 outside $H$, and 
$\nu_{H\backslash G}$ is the invariant measure on 
$H\backslash G$ giving mass 1 to each point, and 
$\dot{k}=Hk$. 
We can rewrite this as 
\vskip.8ex

{\bf Theorem 3.2.} (\hspace{-0.1mm}\cite[\S 3.1]{THirai}) {\it 
The character of\, 
$$
\Pi=\Pi(\pi^0,\pi^1)={\rm Ind}_H^G\pi,\;\;\pi=\pi^0\boxdot\pi^1
$$
 is expressed as  
\begin{gather*}
\chi_\Pi(g) 
 = |H\backslash G|\int_G\chi_\pi(kgk^{-1})\,d\mu_G(k),
\end{gather*}
with the normalized Haar measure $\mu_G$ on $G$, where 
 for $(u,w)\in H=U\rtimes W([\rho])$, \;
$$
\chi_\pi\big((u,w)\big)=\chi_{\pi^0}\big((u,w')\big)\,
\chi_{\pi^1}(w'),
$$ 
with a preimage $w'\in W([\rho])'$ of 
$w$, i.e., $w=\Phi(w')$. 
}%\it 

\vskip.8ex

{\bf 3.3. A complete list of IRs of a semidirect 
product group.} 
For $G=U\rtimes W$, let \,
$$
\big\{\rho_i\;;\;\mbox{\rm IR of $U$},\,i\in I_{U,W}\big\}
$$ 
be a complete set of representatives for the orbit space  
$\widehat{U}/W$, and 
for each $i\in I_{U,W}$, let 
$$
\big\{\pi^1_{i,j}\;;\;j\in J_i\big\}
$$ 
be a complete set of representatives of equivalence classes of  
 IRs of $W([\rho_i])'$ with 
factor set inverse to that of $J_{\rho_i}$. Put \,
\begin{gather*}
\label{}
H_i:=U\rtimes W([\rho_i]),\;H'_i:=U\rtimes W([\rho_i])',
\\
\pi_{i,j}:=\pi^0_i\boxdot \pi^1_{i,j}, 
\\
\Pi_{i,j}:=\Pi(\pi^0_i,\pi^1_{i,j})={\rm Ind}^G_{H_i}\pi_{i,j}.
\end{gather*}

 Define a set of IRs of $G$ as  
 \begin{gather*}
\label{2023-11-30-12}
\Omega(G):=\big\{\Pi_{i,j}
=\Pi(\pi^0_i,\pi^1_{i,j})\;;\;i\in I_{U,W},\, j\in J_i\big\}.
\end{gather*}

{\bf Theorem 3.3.} (\hspace{-0.1mm}\cite[Theorems 3.3 and 4.1]{THirai})  
{\it 
The above set\, $\Omega(G)$ of IRs of\, 
$G$ gives a complete set of representatives of the dual\, 
$\widehat{G}$ of\, $G$, that is, it consists of  
 mutually inequivalent IRs, and covers all equivalence 
classes of IRs exactly once.    
}%\it 

\vskip.8ex

For the proof of this result, the character theory plays 
a decisive roll.

{\bf 4. Structures of $G=G_{20}$ and representation group $R(G)$. } 
Let $G=G_{20}$ be 
a non-abelian group of order 18, presented by 
the pair 
\\
\quad{\bf Set of generators:}\quad$\{a,\,b,\,c\}$;
\\
\quad{\bf Set of fundamental relations:} 
\begin{gather}
\label{2023-12-25-11}
\left\{ 
\begin{array}{ll}
a^2=b^2=c^2=1, 
\\
(abc)^2=1, \quad 
(ab)^3=(ac)^3=1.
\end{array}
\right.
\end{gather}
Its Schur multiplier is  
$M(G_{20})={\mathbb Z}_3$. 
The subgroup $A=\langle a\rangle$ acts on $G$ as
 $\varphi(a)g:=aga^{-1}\;(g\in G)$.

We know by \cite{WP} that $G$ is isomorphic to 
the group $G_{18}^{\;\,4}$ with GAPIdentity [18,4].

{\bf 4.1. Another presentation of $G=G_{18}^{\;\,4}$. } 
We looked for another   
presentation of this group, 
 convenient to our method of constructing 
{\it efficient one-step central extension}, and 
we found a good presentation 
 as follows: 
\\
\quad{\bf Set of generators:} \; $\{x_1=ab,\;x_2=ca,\;a\}$\;;
\\
\quad{\bf Set of fundamental relations:} 
\begin{gather}
\label{2023-11-09-2}
\left\{ 
\begin{array}{ll}
a^2=1, \;x_1^{\;3}=x_2^{\;3}=1, 
\\
x_1x_2=x_2x_1,\;\; \varphi(a)x_i=x_i^{\,-1}\;(i=1,2). 
\end{array}
\right.
\end{gather}

{\bf Lemma 4.1.} {\it Structure of the 
group $G=G_{20}\cong G_{18}^{\;\,4}$ 
is given as follows:\; 
put $X_i:=\langle x_i\rangle\;\;(i=1,2)$, and 
denote by $C_k$ the cyclic group of order $k$,
 then,  
\begin{gather}
\label{2023-12-08-1}
G_{20}=(X_1\times X_2)\rtimes_\varphi A \cong 
(C_3\times C_3)\rtimes_\varphi C_2.
\end{gather}

}%\it

{\bf 4.2. Representation group $R(G_{18}^{\;\,4})$. }  
Taking commuting pair $x_1$ and $x_2$, 
we 
construct a one-step efficient central extension 
$H\stackrel{\delta}{\rightarrow} G$ 
coming from commuting relation $[x_1,x_2]=1$. 
Let $\xi_1, \xi_2, z_{12}$ and $w$ 
be elements of $H$ covering $G$ as 
\begin{gather*}
\xi_i\stackrel{\delta}{\to} x_i\,(i=1,2),\;
w\stackrel{\delta}{\to}a,
\;
[\xi_1,\xi_2]\!=\!z_{12}\stackrel{\delta}{\to} [x_1,x_2]\!=\!1.
\end{gather*}
\indent
The set of all elements of $H$ is 
\begin{gather*}
\label{2023-11-10-1}
h=
z_{12}^{\;\,\beta}\xi_1^{\,\gamma_1}\xi_2^{\,\gamma_2}w^\sigma,
\\
\;\;0\le \beta,\gamma_1,\gamma_2\le 2\;({\rm mod}\;3),\;
0\le \sigma\le 1\;({\rm mod}\;2). 
\end{gather*}  

The set of fundamental relations is 
\begin{gather}
\label{}
\left\{
\begin{array}{ll}
\xi_1^{\;3}=\xi_2^{\;3}=1,\;z_{12}^{\;\,3}=1,\;w^2=1,
\\
z_{12}\;\;{\rm central},\;\;\xi_2\xi_1=z_{12}^{\;-1}\xi_1\xi_2,
\\
\varphi(w)\xi_i=\xi_i^{\,-1}\;(i=1,2). 
\end{array}
\right.
\end{gather}

{\bf Theorem 4.2.} {\it The cyclic subgroup 
$Z_{12}:=\langle z_{12}\rangle$ is the center $Z(H)$ of $H$, 
and isomorphic to Schur multiplier 
$M(G)$ of $G=G_{20}\cong G_{18}^{\;\,4}$.  
The group $H$ 
is a representation group, denoted by $R(G)$. 

Structure of $R(G)$ is given 
with $\Xi_i:=\langle \xi_i\rangle\;\;(i=1,2)$  as 
follows and its GAPIdentity is {\rm [54,8]}: 
$$
R(G)=\big[(\Xi_1\times Z_{12})\rtimes \Xi_2\big]\rtimes W 
\cong [(C_3\times C_3)\rtimes C_3]\rtimes C_2.
$$
}%\it

{\bf Corollary 4.3.} {\it Spin type of an irreducible 
 projective (or spin) representation of $G\cong G_{18}^{\;\,4}$ is 
given as follows: with $\omega:=\exp(2\pi i/3),\,i=\sqrt{-1}$,
$$
\chi_\varepsilon(z_{12})=\omega^\varepsilon, \;\varepsilon =0, \pm 1\;\;\;(
z_{12}\in Z_{12}\cong M(G)). 
$$   

}%\it

{\bf Remark 4.4.} In the list in \cite{WP}, it is noted that 
the group  $G_{54}^{\;\,5}$, with GAPIdentity $[54,5]$, 
has also 
the same 
structure 
$((C_3\times C_3)\rtimes_\varphi C_3)\rtimes_\varphi C_2$. 
Later in the last section, we give more detailed data 
on this point, and show the dual 
(the set of equivalence classes 
of IRs) of $G_{54}^{\;\,5}$.

{\bf 5. Irreducible representations (=IRs) of $G$ itself. } 
Put $U:=X_1\times X_2$ and $W:=A$, 
then $G$ is semidirect product as 
$G=U\rtimes W.$ Since $U$ is abelian, we can 
 apply Mackey's theory here. First, 
the dual $\widehat{U}$ consists of linear characters 
$\rho_{m_1,m_2}$ as follows: 

for 
$(m_1,m_2),\;m_1,m_2\in \{-1,0,1\}\;({\rm mod}\;3)$,  
\begin{gather}
\label{2023-12-09-1}
\rho_{m_1,m_2}(x_k) =\omega^{m_k} \;(k=1,2). 
\end{gather}
The subgroup $W=\langle w\rangle,\;w:=a,$ acts on $\widehat{U}$ as 
\begin{gather}
\label{}
{}^w\rho_{m_1,m_2}\cong \rho_{-m_1,-m_2}, 
\end{gather}
and $\widehat{U}/W$ has two kinds of orbits with parameters   
\\
\;
$(1^*)$ one-point orbit:\;\; $m_1=m_2=0$\,;\quad
\\
\;
$(2^*)$ two-points orbits:\; $m_1\ne 0$ or $m_2\ne 0$. 
\\
Correspondingly the stationary subgroups are   
$$
W([\rho_{0,0}])=W \; \text{and} \; W([\rho_{m_1,m_2}])=\{1\}.
$$ 
\vskip.8ex

\noindent
{\bf Case (1)} \quad \;$\Pi_{0,0;k}:={\bf 1}_U\cdot ({\rm sgn}_W)^k\;(k=0,1)$, \;
\medskip
\\
\quad where ${\bf 1}_U$ is the trivial representation of $U$, and 
\\
\quad
${\rm sgn}_W$ is the sign character of $W\cong \{\pm 1\}$. 
\vskip.8ex

\noindent
{\bf Case (2)} \;Take a representative $\rho=\rho_{m_1,m_2}$  
 and its stationary subgroup $W([\rho])=\{1\}$. Then, 
$$
\Pi_{m_1,m_2}:={\rm Ind}^G_{U\rtimes W([\rho])} \rho_{m_1,m_2}
$$
is an IR looked for. This is realised by using a section 
of $S$ 
 for $U\backslash G \cong S\subset G$ as follows.

Space $V(\Pi),\,\Pi=\Pi_{m_1,m_2}$, consists of $V(\rho)$-valued 
function ${\boldsymbol f}(s),\,s\in S,$ and operation by 
$g_0\in G$ is given by using decomposition  as 
\begin{gather*}
\label{}
sg_0 =u(s,g_0)\,s\overline{g_0},\;\;u(s,g_0)\in U,\;\;
s\overline{g_0}\in S,
\\
\big(\Pi(g_0){\boldsymbol f}\big)(s):=
\rho\big(u(s,g_0)\big)\big({\boldsymbol f}(s\overline{g_0})\big),
\end{gather*} 
where $s\to s\overline{g_0}$ expresses
 the $G$-action on the homogeneous space 
$U\backslash G$. 

{\bf 5.1. Matrix realization of $\Pi_{m_1,m_2}$. }  
Now take $S=W$. Then, with transposed row vectors,
$$
V(\Pi)\ni 
{\boldsymbol f}\to {}^t\big({\boldsymbol f}(1), 
{\boldsymbol f}(w)\big)
\in {\boldsymbol C}^2=:V_2
$$
is an isomorphism. On the space $V_2$, the operators 
$\Pi(u),\,u\in U,$ 
and $\Pi(w)$ are expressed by 
matrices as 
\begin{gather*}
\label{}
\Pi_{m_1,m_2}(u)=
\begin{bmatrix} \rho_{m_1,m_2}(u) &0 \\
0& \rho_{-m_1,-m_2}(u)\end{bmatrix}\!,
\\
\Pi_{m_1,m_2}(w)=
\begin{bmatrix} 0\,&\,1 
\\
1\,&\,0 \end{bmatrix}\!.
\end{gather*}

{\bf Theorem 5.1.} {\it A complete list 
of representatives of the dual of $G\cong G_{18}^{\;4}$ 
is given as follows: 
\\
\;\;{\rm $(1^*)$}\;
two of one-dimensional IRs, 
$\Pi_{0,0;k}\;(k=0,1)$ {\rm ;}
\\
\;\;{\rm $(2^*)$}\;
four of two-dimensional IRs, \;$\Pi_{m_1,m_2}$ with 
\\
\qquad
\;\;$(m_1,m_2)\in\{(1,1), (1,0), (1,-1), (0,1)\}.$

}%\it 
\vskip.8ex

{\bf 5.2. Characters of IRs $\Pi=\Pi_{m_1,m_2}$. } 
We list up trace characters of 2-dimensional IRs 
$\Pi=\Pi_{m_1,m_2}$. The elements of $G$ is expressed 
as 
\begin{gather}
\label{}
\quad
g=x_1^{\,\beta_1}x_2^{\,\beta_2}w^\sigma\;\;
(0\le \beta_1,\beta_2\le 2,\;\sigma=0,1).
\end{gather}

{\bf Theorem 5.2.} {\it The character 
$\chi_\Pi(g)$ is given as  
\begin{gather*}
\label{}
\chi_\Pi(x_1^{\,\beta_1}x_2^{\,\beta_2}w^\sigma)
=
\left\{
\begin{array}{ll}
\omega^{\beta_1m_1+\beta_2m_2}
\\
\qquad +\omega^{-\beta_1m_1-\beta_2m_2} &(\sigma=0), 
\\[.5ex]
\qquad\;\;0&(\sigma=1).
\end{array}
\right.
\end{gather*}
}%\it

{\bf 6. Construction of spin IRs of $G=G_{18}^{\;4}.$ } 
Take the representation group $H:=R(G)$ of $G$, and   
utilize the canonical isomorphism $M(G)\cong Z(H)\cong Z_{12}$.
The set of elements of $H$ is 
\begin{gather*}
\begin{array}{ll}
\qquad\;\; z_{12}^{\;\alpha}\,\xi_1^{\,\beta_1}\xi_2^{\,\beta_2}w^\sigma
\;\; 
(0\le \alpha,\beta_1,\beta_2\le 2,\;\sigma=0, 1).
\end{array}
\end{gather*}

Consider spin IRs $\Pi$ with non-trivial spin type 
\begin{gather}
\label{2023-12-20-1}
\chi_\varepsilon(z_{12})=\omega^\varepsilon\;\;{\rm with}\;\; 
 \;\varepsilon=\pm 1. 
\end{gather}
To construct $\Pi$, we apply the method explained in \S 3 
in four steps and realize it as $\Pi(\pi^0,\pi^1)$ in 
(\ref{2023-11-30-11}).
\vskip.8ex

{\bf 6.1. Subgroup $H_0$ of $H$ and their IRs. } 
To begin with, first 
put 
$U_0=Z_{12}\times \Xi_1$ and take a semidirect 
product subgroup 
$$
H_0=U_0\rtimes \Xi_2\subset H=R(G). 
$$
Then,  $H=H_0\rtimes W$ and we consider two step 
semidirect product operation to arrive until $H=R(G)$. 

First, to construct IRs of $H_0$, we apply 
Mackey's method. The action of $\xi_2\in \Xi_2$ on 
$z_{12},\xi_1\in U_0$ is given as 
$$
\varphi(\xi_2)z_{12}=z_{12},\quad
\varphi(\xi_2)\xi_1=z_{12}^{\;-1}\xi_1.
$$
The dual of\, $U_0$ consists of 
 one-dimensional characters given as follows: 
for $(\varepsilon,m_1),\,\varepsilon, m_1\in\{-1,0,1\}$,
\begin{gather}
\label{2023-12-21-1}
\rho_{\varepsilon,m_1}(z_{12})=\omega^\varepsilon,
\quad\rho_{\varepsilon,m_1}(\xi_1)=\omega^{m_1}.
\end{gather}
Actions of $\xi_2$ on them are 
\begin{gather*}
\varphi(\xi_2)^{-1}z_{12}=z_{12},\;
\varphi(\xi_2)^{-1}\xi_1=z_{12}\xi_1,\;  
\\
\big({}^{\xi_2}\!\rho_{\varepsilon,m_1}\big)(\xi_1)=\omega^{\varepsilon+m_1}
=\rho_{\varepsilon,\varepsilon+m_1}(\xi_1),\;
\\
\big({}^{\xi_2}\!\rho_{\varepsilon,m_1}\big)(z_{12})=\omega^\varepsilon
=\rho_{\varepsilon,\varepsilon+m_1}(z_{12}).\quad
\end{gather*}

{\bf Lemma 6.1.}\; {\it The orbit space 
$\widehat{U_0}/\Xi_2$ consists of 
two kinds of orbits as  
\\
\;\;
{\rm $(3^*)$} three of one point orbits 
$\{\rho_{0,m_1}\}$ with $\varepsilon=0$ {\rm ;} 
\\
\;\;
{\rm $(4^*)$} two of three points 
orbits with $\varepsilon=\pm 1$ {\rm :} 

\qquad ${\cal O}(\varepsilon):=\{\rho_{\varepsilon,m_1};\;m_1=-1,0,1\}$.   

}%\it
\vskip.8ex

Orbits in $(3^*)$ relate with non-spin IRs of $G$, 
since they will 
produce representations with trivial spin type 
$\chi_0(z_{12})=1$. 
So we pick up 
orbits in $(4^*)$.
\vskip.8ex

{\bf 6.2. Induced representation from  $U_0\rtimes \Xi_2([\rho_{\varepsilon,0}])$ to $H_0$. } 
Take $\rho_{\varepsilon,0}$ as a representative of orbits 
${\cal O}(\varepsilon)$, then as is seen from the above, 
the stationary subgroup $\Xi_2([\rho_{\varepsilon,0}])$ in $\Xi_2$ 
is trivial, i.e., 
$\Xi_2([\rho_{\varepsilon,0}])=\{1\}\subset \Xi_2.$
 Therefore, to get IRs of the group $H_0=U_0\rtimes \Xi_2$, 
we should take induced representation from 
$U_0\rtimes \{1\}$ up to $H_0=U_0\rtimes \Xi_2$ as 
$$
\;P_\varepsilon:=
{\rm Ind}^{U_0\rtimes \Xi_2}_{U_0\rtimes \{1\}}
\big(\rho_{\varepsilon,0}\rtimes {\bf 1}_{\{1\}}\big)\;\;{\rm for}\;\,
H_0=U_0\rtimes \Xi_2.
$$

We realize this induced representation on the space of functions 
$f(h)$ of $h\in H_0$ satisfying 
\begin{gather}
\label{2023-12-28-11}
f(uh)=\rho_{\varepsilon,0}(u)f(h)\;\;(u\in U_0,\,h\in H_0),
\end{gather}
with right translation by $h_0\in H_0$ given  as 
$$
f(h)\stackrel{h_0}{\longrightarrow} f(hh_0).
$$

Take a section $S\subset H_0$ for the quotient map 
$H_0\to U_0\backslash H_0$, and restriction map 
$R: f\to \big(f(s)\big)_{s\in S}$ onto the space 
${\cal F}(S)$ of functions on $S$. Then $P_\varepsilon$ 
is expressed as follows: for $h_0\in H_0$,
\begin{gather*}
\;Q_\varepsilon(h_0){:=}R{\cdot} P_\varepsilon(h_0){\cdot}R^{-1}
:
\big(f(s)\big)_{s\in S}\to\big(f(sh_0)\big)_{s\in S}.
\end{gather*}
Choose $S=\Xi_2$ and express 
${}^t\big(f(\xi_2^{\,-1}),f(1),f(\xi_2)\big)$ 
as a row vector (by transposition), then the matrix expression of 
$Q_\varepsilon$ is given as follows. 
\vskip.8ex

{\bf Lemma 6.2.} {\it IRs $Q_\varepsilon$ is given by  
 the following formula: with $I_3$ the identity matrix of order 3,   
\begin{gather*}
\label{}
Q_\varepsilon(z_{12})=\omega^\varepsilon I_3, \quad
Q_\varepsilon(\xi_2)=
{\small 
\begin{bmatrix}
0\,&1\,&0
\\ 
0\,&0\,&1
\\
1\,&0\,&0
\end{bmatrix},
}%\small
\\
Q_\varepsilon(\xi_1)=
\mbox{\rm ${\small 
\begin{bmatrix}
\omega^\varepsilon&0&0
\\
0&\;1\;&0
\\
0&0&\omega^{-\varepsilon}
\end{bmatrix}.
}%\small
$}
\qquad
\end{gather*}
}%\it 

{\it Proof.}\; For $h_0=\xi_2$, we have transformation 
\begin{gather*}
{}^t\big(f(\xi_2^{\,-1}),f(1),f(\xi_2)\big)\to  
{}^t\big(f(1),f(\xi_2),f(\xi_2^{\,-1})\big),
\end{gather*} 
and similarly for $h_0=\xi_1$.
\hfill
$\Box$

\vskip.8ex

{\bf 6.3. From $H_0$ to $H=H_0\rtimes W=R(G)$. } 
Now we apply classical method explained in \S 3, and 
follow from Step 1 to Step 4 in \S 3.1. 

For Step 1, we study action of $w\in W$ to spin IR 
$\rho=Q_\varepsilon$ of $H_0$. Since 
$$
\varphi(w)z_{12}=z_{12},\;\varphi(w)\xi_j=\xi_j^{\,-1}\;(j=1,2),
$$
we have ${}^w\rho(z_{12})=\rho(z_{12})$ and 
$$
{}^w\rho(\xi_j)=\rho\big(\varphi(w)^{-1}\xi_j\big)=\rho(\xi_j^{\,-1})
$$
and we see that ${}^wQ_\varepsilon$ is equivalent to 
$Q_\varepsilon$ and the stationary subgroup is $W([\rho])=W$.

For Step 2, an  
intertwining operator 
$$
J_\rho(w): \delta_{\xi_2^{\,q}}\to \delta_{\xi_2^{\,-q}}\quad{\rm or}\quad 
J_\rho(w)=
{\small 
\begin{bmatrix}
0\,&0\,&1
\\
0\,&1\,&0
\\
1\,&0\,&0
\end{bmatrix}
}%\small
$$ 
between them is given as
\begin{eqnarray*}
\label{2023-12-23-1}
\rho\big(w(h_0)\big)=J_\rho(w)\,\rho(h_0)\,J_\rho(w)^{-1}\;\;(h_0\in H_0). 
\end{eqnarray*}

For Step 3, put for $(h_0, w^\sigma)\in H_0\rtimes W=H$,  
\begin{eqnarray}
\label{2023-12-24-1}
\nonumber
\pi^0\big((h_0,w^\sigma)\big):=\rho(h_0)\cdot J_\rho(w)^\sigma
\;\;\;(\sigma=0,1), 
\end{eqnarray}
then $\pi^0=\rho\cdot J_\rho$, with $\rho=Q_\varepsilon$,  
is an IR of $H$. 

For Step 4,   
take an IR 
$$
\pi^1=({\rm sgn}_W)^k\quad (k=0,1)
$$ 
of $W([\rho])=W$ 
and consider it as a representation of $H$ 
through the homomorphism 
$
H \to\,H/H_0\cong W([\rho])=W,  
$
and consider inner tensor product 
$$
\pi:=\pi^0\boxdot \pi^1
$$ 
as a representation of $H$. 
Since the semidirect product 
$H_0\rtimes W([\rho])$ is just equal to $H=H_0\rtimes W$, 
the inducing process  
${\rm Ind}^{H_0\rtimes W}_{H_0\rtimes W([\rho])}\pi$ for $\pi$ 
is trivial, 
and so we get in this special case 
\begin{gather}
\label{2023-12-24-2}
\Pi(\pi^0,\pi^1)=\pi^0\boxdot \pi^1=:R_{\varepsilon;k}.
\end{gather}

{\bf Proposition 6.3.} {\it Representations 
$R_{\varepsilon;k},\;k=0,1,$ of $R(G)$, with $G=G_{18}^{\;\,4}$, 
is irreducible and of spin type 
 $\chi_\varepsilon$. Their representation matrices are given 
as follows:
\begin{gather}
\label{2023-12-25-1}
\left\{
\begin{array}{rll}
R_{\varepsilon;k}(h_0)&=&Q_\varepsilon(h_0)\quad(h_0\in H_0)
\\   
R_{\varepsilon;k}(w)&=&(-1)^kJ_\rho(w).
\end{array}
\right.
\end{gather}
}%\it

By Theorems 3.1 and 3.2, we obtain  

{\bf Theorem 6.4.} {\it IRs of spin type 
$\chi_\varepsilon$, for $\varepsilon= 1$ or\, 
$\varepsilon=-1$, of 
the group $G=G_{18}^{\;\,4}$ is equivalent to one of 
$R_{\varepsilon;k},\;k=0,1$.

}%\it 

\vskip.8ex

{\bf 6.4. Characters of IRs of $R(G_{18}^{\;\,4})$. } 
All spin IRs $\Pi$ 
of spin type $\chi_\varepsilon,\,\varepsilon=\pm 1,$ 
are 3-dimensional. Therefore their 
trace characters can be calculated 
by explicit representation matrices $\Pi(h),\,h\in H$, and we obtain 
the following formula. 

All elements of $H=R(G_{18}^{\;\,4})$ are expressed as 
\begin{gather*}
\label{}
h=z_{12}^{\;\alpha}\,\xi_1^{\,\beta_1}\xi_2^{\,\beta_2}w^\sigma
\;(
0\le \alpha, \beta_1,\beta_2\le 2,\;\;\sigma=0,1),
\end{gather*} 
and the matrix expression of $R_{\varepsilon;k}(h)$ is 
\begin{gather*}
\label{}
\omega^{\alpha\varepsilon}(-1)^{\sigma k}
{\small 
\begin{bmatrix}
\omega^\varepsilon&0&0
\\0&1&0
\\0&0&\omega^{-\varepsilon}
\end{bmatrix}^{\! \beta_1}\!\!
\begin{bmatrix}
0\,&1\,&0\\
0\,&0\,&1\\
1\,&0\,&0
\end{bmatrix}^{\! \beta_2}\!\!
\begin{bmatrix}
0\,&0\,&1\\
0\,&1\,&0\\
1\,&0\,&0
\end{bmatrix}^{\! \sigma}\!\!\!.
}%\small 
\end{gather*}

{\bf Theorem 6.5.} {\it Character $\chi_\Pi$ for spin IR $\Pi=
R_{\varepsilon;k}$ of spin type $\chi_\varepsilon$ is given for 
$h=z_{12}^{\;\,\alpha}\xi_1^{\,\beta_1}\xi_2^{\,\beta_2}w^\sigma$
as follows:\quad 
 In case $\sigma=0$, 
\\[1ex]
\hspace*{8ex}
$
\chi_\Pi(h)=
\left\{ 
\begin{array}{cl}
3\omega^{\alpha\varepsilon}&\quad  \mbox{\rm 
if $\beta_1=\beta_2=0$}, 
\\
0 &\quad\mbox{\rm otherwise;} 
\end{array}
\right.
$
\\[1ex]
in case $\sigma=1$, 

$
\chi_\Pi(h)=
\left\{ 
\begin{array}{ll} \omega^{\alpha\varepsilon}(-1)^k \quad
& \mbox{\rm if $\beta_2=0$},
\\
\omega^{\alpha\varepsilon}(-1)^k\omega^{-\beta_1\varepsilon}
\quad& \mbox{\rm if $\beta_2=1$},
\\
\omega^{\alpha\varepsilon}(-1)^k\omega^{\beta_1\varepsilon}
& \mbox{\rm if $\beta_2=2$}.
\end{array}
\right.
$
}

\vskip.8ex

{\bf 7. Normalized form of ${\boldsymbol A}=\Pi(a), {\boldsymbol B} =\Pi(b),{\boldsymbol C} =\Pi(c)$. } 
The group $G=G_{18}^{\;\,4}$ is defined firstly 
with {\rm \bf set of generators} $\{a,b,c\}$ and 
{\rm \bf set of fundamental relations} (\ref{2023-12-25-11}) in \S 4. 
Then we introduced a new set of generators as 
$\{x_1=ab,\,x_2=ca,\, a\}$.  
Now, for a spin IR $\Pi$ of $G$, put 
\;$
{\boldsymbol A}=\Pi(a),\; 
{\boldsymbol B} =\Pi(b),\;{\boldsymbol C} =\Pi(c).$\; 
In the previous paper \cite[Theorem 4.1]{HMTY}, it is proved 
that, under certain isomorphisms, the set of 
these matrices can be 
 normalized as follows:
\begin{gather*}
{\boldsymbol A}^2= {\boldsymbol B}^2= 
{\boldsymbol C}^2= I,\; 
({\boldsymbol A}{\boldsymbol B}{\boldsymbol C})^2=\tau I, 
\\
({\boldsymbol A}{\boldsymbol B})^3=
({\boldsymbol A}{\boldsymbol C})^3= I,
\end{gather*}
where $\tau,\,\tau^3=1,$ is a multiplicative constant 
corresponding to 
the spin type of $\Pi$ (hence $M(G)\cong C_3$).
  
Here we list up these 
matrices ${\boldsymbol A},\, 
{\boldsymbol B},\,{\boldsymbol C}$, in this order, 
for every purely-spin IRs $\Pi$ of $G$, 
calculating
  
${\boldsymbol A}=\Pi(w),\; 
{\boldsymbol B} =\Pi(w^{-1}\xi_1),\;
{\boldsymbol C} =\Pi(\xi_2w^{-1})$.
\vskip.5ex

\noindent 
For 1-dim. IR $\Pi_{0,0;k}$: \;$\tau=1$,

\hspace*{3ex} ${\boldsymbol A}=(-1)^k,\;\; 
{\boldsymbol B} =(-1)^k,\;\;
{\boldsymbol C} =(-1)^k$\,;

\noindent 
For 2-dim. IR $\Pi_{m_1,m_2}$: \;$\tau=1$,
{\small 
\begin{gather*}
\quad\;
\begin{bmatrix}
0\;&\;1 
\\
1\;&\;0
\end{bmatrix},\;\;
\begin{bmatrix}
0&\omega^{-m_1}\\
\omega^{m_1}&0 
\end{bmatrix},\;\; 
\begin{bmatrix}
0&\omega^{m_2} 
\\
\omega^{-m_2}&0
\end{bmatrix};
\end{gather*}
}%\small

\noindent
For 3-dim. IRs $R_{\varepsilon;k}$: \;$\tau=\omega^\varepsilon,
\;\varepsilon= 1$ or $-1$,
{\small 
\begin{gather*}
(-1)^k\!
\begin{bmatrix}
0 \! &0 \! &1
\\
0 \!  & 1 \!  &0
\\
1 \! &0 \! &0
\end{bmatrix}\!,\;
(-1)^k\!
\begin{bmatrix}
0 \!  &0 \! &\omega^{-\varepsilon}
\\
0 \!  &1 \! &0
\\
\omega^\varepsilon \! &0 \!  &0
\end{bmatrix}\!,\;
(-1)^k\!
\begin{bmatrix}
0 \!  &1 \! &0
\\
1\!  &0 \! &0
\\
0 \! &0 \! &1
\end{bmatrix}\!.
\end{gather*}
}%\small 

Out of explicit matrix calculations, we see that 
\begin{gather*}
\label{}
{\boldsymbol A}{\boldsymbol B}{\boldsymbol C}
=\Pi(w)\Pi(w^{-1}\xi_1)\Pi(\xi_2w^{-1})
=\Pi(\xi_1\xi_2w^{-1}),
\\
\therefore\quad
({\boldsymbol A}{\boldsymbol B}{\boldsymbol C})^2
=\Pi(\xi_1\xi_2w^{-1}\xi_1\xi_2w^{-1})
\\
=\Pi(\xi_1\xi_2)\Pi(w^{-1}\xi_1\xi_2w)
=\Pi(\xi_1\xi_2\xi_1^{\,-1}\xi_2^{\,-1})
\\
=\Pi(z_{12})=\omega^\varepsilon I_3.
\end{gather*}

{\bf 8. Study on the dual of the group $G_{54}^{\;\,5}$. } 
 In the list in \cite{WP}, it is noted that 
the structure of 
the group  $G_{54}^{\;\,5}$, with GAPIdentity $[54,5]$, 
is expressed as 
$((C_3\times C_3)\rtimes_\varphi C_3)\rtimes_\varphi C_2$, 
the same as for the group $G_{54}^{\;\,8}$. 
So, it would be better or necessary to study 
  group structure of the former and its dual 
 (the set of equivalence classes 
of IRs), to compare with the dual of the latter. 
Rewriting data from GAP system we have the following 
presentation (cf. \cite[p.4]{Tsuri0}): 

\noindent
\quad{\bf Set of generators:} \quad $\{h_1,\,h_2,\,h_3,\,h_4\}$;
\\
\quad
{\bf Set of fundamental relations:}   
\begin{gather*}
\left\{
\begin{array}{ll}
h_1^{\,2}=1,\;h_2^{\,3}=h_3^{\,3}=h_4^{\,3}=1,
\\
\varphi(h_1)h_2=h_2,\;\varphi(h_2)h_4=h_4,\;\varphi(h_3)h_4=h_4,
\\
\varphi(h_1)h_3=h_3^{-1}\!,\,\varphi(h_1)h_4=h_4^{-1}\!,\,
\varphi(h_2)h_3=h_3h_4, 
\end{array}
\right.
\end{gather*}
where $\varphi(h)x=hxh^{-1}$. 
Denote by $H_j:=\langle h_j\rangle \;(1\le j\le 4)$ the cyclic group 
generated by $h_j$, and 
$$
W{:=}H_1{\times} H_2\cong C_2{\times} C_3,
\;\,U{:=}H_3{\times} H_4\cong C_3{\times} C_3, 
$$
then 
$U$ is an abelian normal subgroup 
and $W$ acts on $U$ as shown 
in the fundamental relations. So the structure of 
$G_{54}^{\;\,5}$ is described 
more clearly as 
$$
G_{54}^{\;\,5}\cong (C_3\times C_3)\rtimes (C_2\times C_3).
$$

To construct a complete set of representatives (IRs) 
of the dual of 
this group, we can also apply Mackey's induced 
representations, 
since $U$ is abelian.

\noindent
{\large \bf e-mail addresses:}

\medskip

\hspace{-2mm}Tatsuya Tsurii: t3tsuri23@rsch.tuis.ac.jp

\hspace{-2mm}Satoe Yamanaka: yamanaka@libe.nara-k.ac.jp

\hspace{-2mm}Itsumi Mikami: kojirou@kcn.ne.jp

\hspace{-2mm}Takeshi Hirai: hirai.takeshi.24e@st.kyoto-u.ac.jp

\end{multicols*}

\end{document}